# Causality and Statistical Learning[1]

Andrew Gelman[2]

24 Feb 2010

In social science we are sometimes in the position of studying descriptive questions (for example: In what places do working-class whites vote for Republicans? In what eras has social mobility been higher in the United States than in Europe? In what social settings are different sorts of people more likely to act strategically?). Answering descriptive questions is not easy and involves issues of data collection, data analysis, and measurement (how should one define concepts such as "working class whites," "social mobility," and "strategic"), but is uncontroversial from a statistical standpoint.

All becomes more difficult when we shift our focus from What to What-if and Why.

**Thinking about causal inference**

Consider two broad classes of inferential questions:

1. *Forward causal inference*. What might happen if we do X? What are the effects of smoking on health, the effects of schooling on knowledge, the effect of campaigns on election outcomes, and so forth?

2. *Reverse causal inference*. What causes Y? Why do more attractive people earn more money, why do many poor people vote for Republicans and rich people vote for Democrats, why did the economy collapse?

In forward reasoning, the potential treatments under study are chosen ahead of time, whereas, in reverse reasoning, the research goal is to find and assess the importance of the causes. The distinction between forward and reverse reasoning (also called "the effects of causes" and the "causes of effects") was made by Mill (1843). Forward causation is a pretty clearly-defined problem, and there is a consensus that it can be modeled using the counterfactual or potential-outcome notation associated with Neyman (1923) and Rubin (1974) and expressed using graphical models by Pearl (2009): the causal effect of a treatment T on an outcome Y for an individual person (say), is a comparison between the value of Y that would've been observed had the person followed the treatment, versus the value that would've been observed under the control; in many contexts, the treatment effect for person i is defined as the difference, $Y_i(T=1)$ -

[2] Department of Statistics and Department of Political Science, Columbia University, New York, gelman@stat.columbia.edu, http://www.stat.columbia.edu/gelman/  We thank Kevin Quinn, Jasjeet Sekhon, Judea Pearl, Rajeev Dehejia, Michael Sobel, Jennifer Hill, Paul Rosenbaum , Sander Greenland, Guido Imbens, Donald Rubin, Steven Sloman, Peter Hoff, Thomas Richardson, Niall Bolger, and Christopher Winship for helpful comments and the Institute for Educational Sciences for partial support of this work.

$Y_i(T=0)$.  Many common techniques, such as differences in differences, linear regression, and instrumental variables, can be viewed as estimated average causal effects under this definition.

In the social sciences, where it is generally not possible to try more than one treatment on the same unit (and, even when this is possible, there is the possibility of contamination from past exposure and changes in the unit or the treatment over time), questions of forward causation are most directly studied using randomization or so-called natural experiments (see Angrist and Pischke, 2008, for discussion and many examples). In some settings, crossover designs can be used to estimate individual causal effects, if one accepts certain assumptions about treatment effects being bounded in time.  Heckman (2006), pointing to the difficulty of generalizing from experimental to real-world settings, argues that randomization is not any sort of "gold standard" of causal inference, but this is a minority position:  I believe that most social scientists and policy analysts would be thrilled to have randomized experiments for their forward-causal questions, even while recognizing that subject-matter models are needed to make useful inferences from any experimental or observational study.

Reverse causal inference is another story.  As has long been realized, the effects of action X flow naturally forward in time, while the causes of outcome Y cannot be so clearly traced backward.  Did the North Vietnamese win the American War because of the Tet Offensive, or because of American public opinion, or because of the skills of General Giap, or because of the political skills of Ho Chi Minh, or because of the conflicted motivations of Henry Kissinger, or because of Vietnam's rough terrain, or . . .?  To ask such a question is to reveal the impossibility of answering it.  On the other hand, questions such as "Why do whites do better than blacks in school?", while difficult, do not seem inherently unanswerable or meaningless.

We can have an idea of going backward in the causal chain, accounting for more and more factors until the difference under study disappears—that is, is "explained" by the causal predictors.  Such an activity can be tricky—hence the motivation for statistical procedures for studying causal paths—and ultimately is often formulated in terms of forward causal questions: causal effects that add up to explaining the Why question that was ultimately asked.  Reverse causal questions are often more interesting and motivate much, perhaps most, social science research; forward causal research is more limited and less generalizable but is more doable.  So we all end up going back and forth on this.

We see three difficult problems in causal inference:

- Studying problems of forward causation with observational data or experiments with missing data (the traditional focus of causal inference in the statistics and biostatistics literature); recall that missingness is inherent in the counterfactual definition of causal effects.

- Generalizing from experiments or quasi-expeirments to realistic settings. This is the problem of external validity (Campbell and Stanley, 1966), which Heckman (2006) has emphasized.  The difficulty is that experiments are often a poor method for understanding how effects may vary across contexts, including from the lab to the real world, and across populations when treatment effects vary.

- Studying reverse causation (learning about "Why" questions) from observational data using multivariate analysis. Here the many threats to a model's validity may be so great as to make confidence in any estimates unattainable..

Along with any of these analyses comes the problem of assessing the plausibility of whatever assumptions are being used to identify causal effects in the particular problem at hand, and the assessment of the sensitivity of causal inferences to these assumptions.

It was with this in mind that I read the three books under review, each of which is well-written, thought-provoking, and advances the science of causal reasoning. Morgan and Winship focus on the second item above (their running examples include the effects of family background on educational attainment, the effects of education on earnings, and the effects of voting technology on voting), and Sloman also talks a bit about the third item: learning about causal relationships from multivariate observational data. Fernbach, Darlow, and Sloman (2010) discuss the role of forward and reverse causal reasoning in everyday judgment. Pearl presents a framework in which the potential-outcome model for forward causation can be used to study reverse causation as well. (The three books pretty much ignore the first item in my list above, perhaps under the assumption that, once the basic causal questions are identified, it shouldn't be too hard to fit a model and use it to generalize to other populations.)

In this review, I address some of the central debates in causal inference but I do not attempt anything like a serious literature review. I encourage readers to follow the trails of references in the cited here to get a fuller sense of the development of causal inference ideas in the social and cognitive sciences.

The statistical literature on causal inference (which has been influenced by work in economics, epidemiology, psychology, philosophy, and computer science) has featured contentious debate, and some of that debate has found its way into the books under review here. I recommend starting with Morgan and Winship's *Counterfactuals and Causal Inference*, which is written by two sociologists and provides an agnostic overview of the field. They discuss regression estimates of causal effects and the assumptions needed for these estimates to make sense. Along with this, I recommend the discussion of Sobel (2008), who discusses the assumptions required to make inference for mediation—that is, quantitative estimates of the decomposition of a causal effect across different paths—and Rosenbaum's (2010) book on the design of studies for causal inference. In social science, it can be important not just to have one substantive theory, but to have two theories which the study will contrast, and to have a design that can evaluate them.

Pearl's book is the second edition of his influential work that presents causal reasoning as mathematical operations on a system of non-parametric equations, often encoded in directed graphical models. Pearl's formulation is, at the very minimum, a helpful way to understand structural equations and identifiability (as explained by Morgan and Winship in language that may be more comfortable to a statistically-trained social scientist), and this is important. Without a good formal structure, we are all at risk of getting confused over causal inference; for instance, Rubin, 2005, relates an example where the legendary R. A. Fisher made the basic (but apparently not trivial) mistake of adjusting for an intermediate outcome when estimating a treatment effect. Even if we are alert to this particular mistake, if we work without a formal structure we are still at risk of making others, for example when deciding whether to adjust for particular pre-treatment variables in the presence of unobserved confounders.

The third book under review is different from the others. Steven Sloman is a cognitive scientist and is just as interested in how people *do* think as in how they *should* think. And people think about causality all the time. Just as Karl Popper (1959) and others have told us, scientific facts are difficult to interpret or even collect except in light of a theory.

In his research, Sloman has studied the ways that people use causal models to understand correlations. For example:

> A fact increases belief in another fact if the two share an explanation. For example, telling people that boxers are not eligible for most health insurance plans increased their judgment of the probability that asbestos removers are not eligible either. The common explanation depends on the causal relations.
> (Increased health risk → Cost to health insurers → Denial of coverage)
> However, telling people that boxers are more likely than average to develop a neurological disorder actually decreased their average probability that asbestos removers are more likely than average to develop a neurological disorder.

This sort of finding demonstrates the importance of explanation in our understanding of the world. Although it does not itself represent statistical learning, the evident importance of causal thinking in everyday (not just scientific or statistical) contexts makes cognitive and information scientists receptive to Pearl's view that causal structure can, under certain conditions, be perceived from correlational data.

**Attitudes toward causal inference**

We can identify a (very rough) ordering of views based on the sorts of assumptions used in causal reasoning. At one extreme is a view expressed by the economists James Heckman (2006), who trusts causal inferences only when they are based on a substantively-motivated model (although this stance can be more permissive to the extent that one is willing to accept such substantive models); see also Deaton (2009). Slightly less conservative is the mainstream view in statistics and labor economics (see Angrist and Pischke, 2008) in which causal inference is allowed from randomized experiments and in certain quasi-random or "natural experimental" settings under strong assumptions about data collection (such as associated with regression discontinuity designs or instrumental variables). In that way, the statisticians and econometricians, even when they use Bayesian inference, are working with the classical statistical paradigm in which identification is driven by randomization (or, more generally, conditional independence) in data collection. Compared to other social scientists, economists tend to focus more on the effects of possible interventions. Moving in the more permissive direction, we find many epidemiologists (for example, Greenland and Robins, 1986, 2009) who, having long experience with procedures such as attributable risk analysis, are more comfortable with causal inference from observational data, feeling that the assumptions required to move from association to causation are reasonable in some scientific settings. Even more permissive

(in this admittedly oversimplified ordering) are social psychologists and others who have been using path analysis and structural equation models for many years to identify the strengths of posited causal relationships from observed covariance matrices (see Duncan, 1966, and Kenny, 1979)—although in recent years there has been a move toward a hybrid approach in which mediation analysis is applied to experimental data to study the effects of intervening variables (Baron and Kenny, 1986). Finally, the extreme position is taken, I believe, by those who believe that a computer should be able to discern causal relationships from observational data, based on the reasonable argument that we, as humans, can do this ourselves in our everyday life with little recourse to experiment.

The books under review stand somewhat outside this spectrum. Morgan and Winship only commit to the conservative statistical view—causal inference from experiments and quasi-experiments—but they present the logic behind the more permissive stances. Pearl, as befits his computer-science background, supports a permissive view, but only in theory; he accepts that any practical causal inference will be only as good as the substantive models that underlie it and he is explicit about the additional assumptions required to translate associations into causal inference. In some ways he is close to the statisticians and econometricians who emphasize that randomization allows more robust inferences for causal inference about population averages, although he is also interested in causal inference when randomization or its equivalents are not available. (Pearl does discuss methods for learning causal structural from observational data alone, but such inference is done under an assumption called faithfulness or stability, which corresponds to a model of conditional independence, which I do not think makes sense in social science settings except when such independence arises by design; see the discussion of "no true zeroes" below.)

I myself hold conflicting views, in a manner similar to that expressed by Sobel (2008). On one hand, I do not see how one can get scientifically-strong causal inferences from observational data alone without strong theory; it seems to me a hopeless task to throw a data matrix at a computer program and hope to learn about causal structure. (More on this below.) On the other hand, I recognize that much of our everyday causal intuition, while not having the full quality of scientific reasoning, is still useful, and it seems a bit of a gap to simply use the label "descriptive" for all inference that is not supported by experiment or very strong theory.

Thus, I am sympathetic to theories such as Pearl's that bridge different notions of causal inference and I applaud research such as Sobel's that, by clarifying implicit assumptions, bridges the gap between the experimental and structural-modeling perspectives, even as I have difficulty seeing exactly how to apply them in my own research.

For example, income, religion, and religious attendance predict vote choice in different ways in different parts of the country, and we have also found variation in attitudes on economic and social issues (Gelman et al., 2009). But I don't know exactly what would be learned by throwing all these variables into a multivariate model. Throwing a data matrix into causal multivariate analysis software (such as that developed from the ideas of Spirtes, Glymour, and Scheines, 2001) is, if nothing else, a form of data summary, and I can well believe that it might be useful to those trained to interpret that data reduction, even if the estimated graphical model and estimated coefficients can't be interpreted directly in the way that the designers of these programs might think.  I am open to the possibility that such an analysis could be useful, but we don't appear to be there yet, and the outputs of structural modeling programs can easily be misinterpreted.

**There are (almost) no true zeroes:  difficulties with the research program of learning causal structure**

We can distinguish between learning within a causal model (that is, inference about parameters characterizing a specified directed graph) and learning causal structure itself (that is, inference about the graph itself).  In social science research, I am extremely skeptical of this second goal.

The difficulty is that, in social science, there are no true zeroes.  For example, religious attendance is associated with attitudes on economic as well as social issues, and both these correlations vary by state.  And it does not interest me, for example, to test a model in which social class affects vote choice through party identification but not along a direct path.  More generally, anything that plausibly *could* have an effect will not have an effect that is exactly zero.  I can respect that some social scientists find it useful to frame their research in terms of conditional independence and the testing of null effects, but I don't generally find this approach helpful—and I certainly don't believe that it is *necessary* to think in terms of conditional independence in order to study causality.  Without structural zeroes, it is impossible to identify graphical structural equation models.

The most common exceptions to this rule, as I see it, are independences from design (as in a designed or natural experiment) or effects that are zero based on a plausible scientific hypothesis (as might arise, for example, in genetics (where genes on different chromosomes might have essentially independent effects, or in a study of ESP).  In such settings I can see the value of testing a null hypothesis of zero effect, either for its own sake or to rule out the possibility of a conditional correlation that is supposed not to be there.

Another sort of exception to the "no zeroes" rule comes from information restriction:  a person's decision should not be affected by knowledge that he or she doesn't have.  For example, a consumer interested in buying apples cares about the total price he pays, not about how much of that goes to the seller and how much goes to the government in the form of taxes.  So the restriction is that the utility depends on prices, not on the share of that going to taxes.  That is the type of restriction that can help identify demand functions in economics.

I realize, however, that my perspective that there are no zeroes (information restrictions aside) is a minority view among social scientists and perhaps among people in general, on the evidence of psychologist Sloman's book.  For example, from chapter 2:  "A good politician will know who is motivated by greed and who is motivated by larger principles in order to discern how to solicit each one's vote when it is needed."  I can well believe that people think in this way but I don't buy it!  Just about everyone is motivated by greed *and* by larger principles!  This sort of discrete thinking doesn't seem to me to be at all realistic about how people behave—although it might very well be a good model about how people characterize others!

In the next chapter, Sloman writes, "No matter how many times A and B occur together, mere co-occurrence cannot reveal *whether* A causes B, or B causes A, or something else causes both."

[italics added]  Again, I am bothered by this sort of discrete thinking.  I will return in a moment with an example, but just to speak generally, if A *could* cause B, and B *could* cause A, then I would think that, yes, they could cause each other.  And if something else *could* cause them both, I imagine that could be happening along with the causation of A on B and of B on A.

Here we're getting into some of the differences between a normative view of science, a descriptive view of science, and a descriptive view of how people perceive the world.  Just as there are limits to what "folk physics" can tell us about the motion of particles, similarly I think we have to be careful about too closely identifying "folk causal inference" from the stuff done by the best social scientists.  To continue the analogy:  it is interesting to study how we develop physical intuitions using commonsense notions of force, energy, momentum, and so on—but it's also important to see where these intuitions fail.  Similarly, ideas of causality are fundamental but that doesn't stop ordinary people and even experts from making basic mistakes.

Now I would like to return to the graphical model approach described by Sloman.  In chapter 5, he discusses an example with three variables:

> If two of the variables are dependent, say, intelligence and socioeconomic status, but conditionally independent given the third variable [beer consumption], then either they are related by one of two chains:
> (Intelligence → Amount of beer consumed → Socioeconomic status)
> (Socio-economic status → Amount of beer consumed → Intelligence)
> or by a fork:
> (Amount of beer consumed → Socioeconomic status)
>                            → Intelligence)
> and then we must use some other means [other than observational data] to decide between these three possibilities.  In some cases, common sense may be sufficient, but we can also, if necessary, run an experiment.  If we intervene and vary the amount of beer consumed and see that we affect intelligence, that implies that the second or third model is possible; the first one is not.  Of course, all this assumes that there aren't other variables mediating between the ones shown that provide alternative explanations of the dependencies.

This makes no sense to me.  I don't see why only one of the three models can be true.  This is a mathematical possibility, but it seems highly implausible to me.  And, in particular, running an experiment that reveals one of these causal effects does *not* rule out the other possible paths.  For example, suppose that Sloman were to perform the above experiment (finding that beer consumption affects intelligence) and then *another* experiment, this time varying intelligence (in some way; the method of doing this can very well determine the causal effect) and finding that it affects the amount of beer consumed.

Beyond this fundamental problem, I have a statistical critique, which is that in social science you won't have these sorts of conditional independencies, except from design or as artifacts of small sample sizes that do not allow us to distinguish small dependencies from zero.

I think I see where Sloman is coming from, from a psychological perspective:  you see these variables that are related to each other, and you want to know which is the cause and which is the effect. But I don't think this is a useful way of understanding the world, just as I don't think it's useful to categorize political players as being motivated either by greed or by larger principles, but not both. Exclusive-or might feel right to us internally, but I don't think it works as science.

One important place where I agree with Sloman (and thus with Pearl and Sprites et al.) is in the emphasis that causal structure cannot in general be learned from observational data alone; they hold the very reasonable position that we can use observational data to rule out possibilities and formulate hypotheses, and then use some sort of intervention or experiment (whether actual or hypothetical) to move further. In this way they connect the observational/experimental division to the hypothesis/deduction formulation that is familiar to us from the work of Popper, Kuhn, and other modern philosophers of science.

The place where I think Sloman is misguided is in his formulation of scientific models in an either/or way, as if, in truth, social variables are linked in simple causal paths, with a scientific goal of figuring out if A causes B or the reverse. I don't know much about intelligence, beer consumption, and socioeconomic status, but I certainly don't see any simple relationships between income, religious attendance, party identification, and voting—and I don't see how a search for such a pattern will advance our understanding, at least given current techniques. I'd rather start with description and then go toward causality following the approach of economists and statisticians by thinking about potential interventions one at a time. I'd love to see Sloman's and Pearl's ideas of the interplay between observational and experimental data developed in a framework that is less strongly tied to the notion of choice among simple causal structures.

**Causality as intervention and causality as system modeling**

The present discussion has highlighted some different ways of thinking about causality. In the world of Pearl, causation is defined in terms of interventions (the "do" operator), and causality is restricted to be unidirectional, forward in time—hence the *directed* graphs. Rubin adds another restriction, that of associating any causal effect with to a particular intervention, an idea that Sekhon (2009), via Cochran (1965), traces to the following question of Dorn (1956): "How would the study be conducted if it were possible to do it by controlled experimentation?"

For example, statistical writing on causality has sometimes been misunderstood to imply that you can't ever estimate the effects of sex and ethnicity, because it is (nearly) impossible to manipulate these on an individual person. The real point, though, is not that such effect can't be studied but rather that there are many versions of such effects. One sort of effect of sex and ethnicity on hiring has been studied using experiments in which interviewers have been asked to evaluate files with and without details on the job candidates. In this case, the effects are identified with this one particular sort of intervention; other, biological, effects of sex and ethnicity could be studied in other ways.

Perhaps because of my own statistical training, I am most comfortable with this view, in which causal inferences are tied to specific (real or hypothesized) treatment. Thus, for example, one might consider the effect of a law requiring adolescents to stay in school another year, or the effect of a financial incentive to stay in school, rather than the more generic effect of schooling.

I find it helpful to focus ideas by thinking of specific treatments, but one could certainly argue the opposite and hold the position that the graphical modeling approach (with its "do" operator that can be applied to any variable in a system) is a more realistic way to understand the world as it is.

Different from both the above approaches is the world of many structural equation models used in social science, where causation is defined in terms of variables in a system (for example, inflation can affect unemployment, unemployment can affect inflation, and both can both affect and be affected by government policy). The system-variable approach can be incorporated into the do-calculus or intervention framework by adding time into the system (inflation at time 1 affects unemployment at time 2, and so forth), but this is not always so easy to do with observational data, and in many ways the goals are different: unidirectional causation in one case and equilibrium modeling in the other.

Sloman's book is interesting partly because it points up the tension between these two views, both of which are present in our casual reasoning about causality. If we take three variables (beer consumption, intelligence, and social status) that are defined in no particular logical or time order, then any causal patterns are possible—and they can all hold. From the causality-as-intervention perspective, it makes sense to study the effects of each of these three variables considered as a treatment, and this will likely require three different experiments or observational studies. It would rely on unrealistic assumptions to think that one observational analysis could get at all three treatment effects. From the system-variable perspective, we can consider the effects of these variables on each other without imagining that one configuration of the arrows is correct and the others wrong.

Similarly, in economics and sociology, researchers have studied the returns to education and the factors that influence students to stay in school. I find it more helpful to understand the causality here in terms of two different sorts of interventions rather than as a single set of simultaneous equations.

**Conclusions**

Casual inference will always be a challenge, partly because our psychological intuitions do not always match how the world works. We like to think of simple causal stories, but in the social world, causes and effects are complex. We—the scientific community—still have not formed a consensus about how best to form consensus on causal questions: there is a general acceptance that experimentation is a gold standard, but it is not at clear how much experimentation should be done, to what extent we can trust inference from observational data, and to what extent experimentation should be more fully incorporated into daily practice (as suggested by some in the "evidence-based medicine" movement).

The most compelling causal studies have (i) a simple structure that you can see through to the data and the phenomenon under study, (ii) no obvious plausible source of major bias, (iii) serious efforts to detect plausible biases, efforts that have come to naught, and (iv) insensitivity to small and moderate biases (see, for example, Greenland, 2005). Two large open questions are, first, how to best achieve these four steps in practice and, second, what sorts of causal claims to make in settings where we are not able to satisfy these conditions. Right now I think the best approach is a combination of the economists' focus on clever designs and identification strategies and the statisticians' ability to build more complicated models to assess what might happen if the strict assumptions fall apart.

The three books under review convey much of the state of the art in causal reasoning. My own view is that for the foreseeable future, those of us who do causal reasoning in the social sciences will have to continue to think hard about identification strategies (as discussed, for example, by Campbell and Stanley, 1966, and Angrist and Pischke, 2008) and to recognize that multiple analyses may be needed to explore complex causal paths. Ideas of discovering an underlying causal structure do not match to how I understand social science.

But there is a diversity of perspectives, both in theory and in practice, on how to move from description to causal inference, and I recommend that social scientists try different approaches as befit their applications. The research contributions of Pearl, and the exposition of Morgan and Winship, should give readers the tools to state the assumptions underlying whatever causal models they choose to use, and the research of Sloman and other cognitive scientists should help us in better understanding how people draw causal conclusions in general, and in the relations between observational and experimental inference inside and outside of science.